\documentclass[12pt,draftcls,onecolumn]{IEEEtran}
\usepackage{fancyhdr}

\def\baselinestretch{1.2}

\usepackage{amsfonts}
\usepackage{mathrsfs}
\usepackage{}
\usepackage{colortbl}
\usepackage{amssymb,amsmath,graphicx}

\usepackage{caption}
\usepackage{subfigure}
\usepackage{float}

\newtheorem{Remark}{\bf Remark}[section]

\newtheorem{Problem}{\bf Problem}[section]

\newenvironment{Proof}{\noindent{\em Proof:\/}}{\hfill $\Box$\par}

\newtheorem{Theorem}{\bf Theorem}[section]

\newtheorem{Lemma}{\bf Lemma}[section]

\newtheorem{Assumption}{\bf Assumption}[section]

\newcommand{\EQ}{\begin{eqnarray}}

\newcommand{\EN}{\end{eqnarray}}

\newcommand{\EQQ}{\begin{eqnarray*}}

\newcommand{\ENN}{\end{eqnarray*}}

\begin{document}

\title{\LARGE
Adaptive Leader-Following Consensus for a Class of Higher-Order Nonlinear Multi-Agent Systems with Directed Switching Networks
}



\author{Wei~Liu~and~Jie~Huang
\thanks{The original version of this paper appeared recently in \cite{LH5}. The main difference of this version from \cite{LH5} is that  Appendix B is added to show the existence of the limit of the function $V (t)$ defined in \eqref{V3} as $t$ tends to infinity.}
\thanks{Wei Liu and Jie Huang are with
 Department of Mechanical and Automation Engineering, The Chinese University of Hong Kong, Shatin, N.T., Hong Kong. Email: wliu@mae.cuhk.edu.hk, jhuang@mae.cuhk.edu.hk.}
\thanks{Corresponding author: Jie Huang.}
}

\maketitle

\thispagestyle{fancy}
\fancyhead{}
\lhead{}
\lfoot{}
\cfoot{}
\rfoot{}
\renewcommand{\headrulewidth}{0pt}
\renewcommand{\footrulewidth}{0pt}

\begin{abstract}
In this paper, we  study the leader-following consensus problem for a class of uncertain nonlinear multi-agent systems under jointly connected directed switching networks.
The uncertainty includes constant unbounded parameters and external disturbances. We first extend the recent result on the adaptive distributed observer from global asymptotical convergence to global exponential convergence. 
 Then,  by integrating the conventional adaptive control technique with the adaptive distributed observer, we present our solution by a distributed adaptive state feedback control law.
Our result is illustrated by the leader-following consensus problem for a group of van der Pol oscillators.
\end{abstract}

\section{Introduction}\label{Introduction}
In the past few years, the cooperative control problems for multi-agent systems have attracted extensive attention due to their wide applications in engineering systems such as sensor networks, robotic teams, satellite clusters, unmanned air vehicle formations and so on. The consensus problem is one of the basic cooperative control problems, whose objective is to design a distributed control law for each agent such that the states (or outputs) of all agents approach the same value. Depending on whether or not a multi-agent system has a leader, the consensus problem can be divided into
two classes: leaderless and leader-following. The leaderless consensus problem aims to make the states (or outputs) of all agents asymptotically synchronize to a same trajectory, while the leader-following consensus problem requires the states (or outputs) of all agents to asymptotically track a desired trajectory which is generated by the leader system.

The consensus problem of linear multi-agent systems has been extensively studied. For example,  the leaderless case was studied in \cite{Olfati1,Ren1,Seo1,Tuna1},  the leader-following case was studied in \cite{Hong1,Hong2,Hu1,Hu3,Ni1},  and  both two cases were studied in \cite{Jadbabaie1,SH3}. In particular, the linear multi-agent system considered in \cite{Hu3} contains some time-varying disturbances, and  the adaptive control technique has been used to deal with these disturbances. Recently, more attention has been paid to the consensus problem of nonlinear multi-agent systems. For example, 
in \cite{LiuK1,Mei1,Song1,Song2,Yu1}, the consensus problem was studied for several classes of nonlinear systems 
satisfying the global Lipschitz condition or the global Lipschitz-like condition.
 In \cite{LH1,LH4,SH4,WangX1}, the leader-following consensus problem was studied  via the output regulation theory and the nonlinear systems considered in \cite{LH1,LH4,SH4,WangX1} contain both disturbance and uncertainty, but the boundary of the uncertainty is known.  In \cite{Hu2},  the authors designed a nonlinear observer-based filter to track a single second-order linear Gaussian target and analyzed the stability of the proposed filter in the sense of mean square. Based on the adaptive control technique, the leader-following consensus problem  was studied  for  first-order  nonlinear multi-agent systems  in \cite{Das1,YuH1}, for second-order nonlinear multi-agent systems in \cite{LH3}, and for multiple uncertain rigid spacecraft systems in \cite{CH2}. In \cite{Das1,ZhangH1},  the neutral networks method was used to study uncertain nonlinear multi-agent systems subject to static networks and the designed control laws can make the tracking errors uniformly ultimately bounded  for initial conditions in some prescribed compact subset.

In this paper, we will further consider the leader-following consensus problem for a class of uncertain nonlinear multi-agent systems. Our paper has the following features. First,
the order of our system is generic and the nonlinearity does not have to satisfy the global Lipschitz-like condition which  excludes some benchmark nonlinear systems such as van der Pol systems, Duffing systems and so on. Thus, the linear control techniques as used in
 \cite{LiuK1,Mei1,Song1,Song2,Yu1} do not apply to our system. Second, our system contains both constant uncertain parameters and external disturbances and  the uncertain parameters can take any constant  value. Thus, the robust control approaches in \cite{LH1,LH4,SH4,WangX1} do not apply to our system either. Third,  our networks satisfy the jointly connected condition, which is the mildest condition on the communication network since it allows the network to be disconnected at any time, and  contains the static network case \cite{Das1,ZhangH1} and the every time connected switching network case \cite{LH1} as special cases.
Finally, compared with \cite{Das1,ZhangH1}, our result is global and the consensus can be achieved exactly.
As a result of these features, the problem is much more general than the existing results and cannot be handled by the techniques in the literatures.
To solve our problem, we have integrated the classical adaptive control technique and the recently developed adaptive distributed observer to obtain a distributed adaptive control law.
We have also  furnished a detailed stability analysis for the closed-loop system.

It should be noted that the leader-following consensus problem for a class of multiple uncertain Euler-Lagrange systems has been studied in \cite{CH1}, where the adaptive distributed observer method has been first proposed. However, the system considered in \cite{CH1} contains only parameter uncertainty but no disturbance,  and the communication network is assumed to be  undirected jointly connected. In this paper, we extend the network from the undirected case to the directed case.

The rest of this paper is organized as follows. In Section \ref{PF}, we present our problem formulation and two assumptions. In Section \ref{Observer}, we introduce some concepts for the adaptive distributed observer and establish a technical lemma. In Section \ref{MR}, we present our main result. In Section \ref{Example},
we provide an example to illustrate our design. Finally, in Section \ref{Conclusion}, we close the paper with some concluding remarks.

{\bf Notation.} For any column vectors $a_i$, $i=1,...,s$, denote $\mbox{col}(a_1,...,a_s)=[a_1^T,...,a_s^T]^T$. $\otimes$ denotes the Kronecker product of matrices. Vector $\mathbf{1}_{N}$ denotes an $N$-dimensional column vector with all elements being $1$.
$\|x\|$ denotes the Euclidean norm of vector $x$. $\|A\|$ denotes the induced norm of matrix $A$ by the Euclidean norm.
 $\lambda_{\max}(A)$ and $\lambda_{\min}(A)$ denote the maximum eigenvalue and the minimum eigenvalue of a symmetric real matrix $A$, respectively.
 We use $\sigma(t)$ to denote a piecewise constant switching signal $\sigma:[0,+\infty)\rightarrow \mathcal{P}=\{1,2,\dots,n_{0}\}$, where $n_{0}$ is a positive integer, and
$\mathcal{P}$ is called a switching index set. We assume that all switching instants $t_0=0 <t_1<t_2,\dots$ satisfy $t_{i+1}-t_{i}\geq \tau_{0}>0$ for some constant $\tau_{0}$ and all $i\geq 0$, where $\tau_{0}$ is called the dwell time.


\section{Problem Formulation}\label{PF}
Consider a class of nonlinear multi-agent systems as follows:
\begin{equation}\label{system1}
\begin{split}
    \dot{x}_{si}& =x_{(s+1)i},~s=1,2\cdots,r-1\\
    \dot{x}_{ri}& = f_{i}^{T}(x_{i},t)\theta_{i}+d_{i}(w)+u_{i},~~ i=1,\cdots,N
\end{split}
\end{equation}
where $x_{i}=\mbox{col}(x_{1i},\cdots,x_{ri})\in\mathbb{R}^{r}$ is the state,  $u_{i}\in \mathbb{R}$ is the input, $f_{i}:\mathbb{R}^{r}\times[0,+\infty)\rightarrow\mathbb{R}^{m}$ is a known function satisfying locally Lipschitz condition with respect to $x_{i}$ uniformly in $t$,  $\theta_{i}\in\mathbb{R}^{m}$ is an unknown constant parameter vector, $d_{i}(w)$ denotes the disturbance with  $d_{i}:\mathbb{R}^{n_{w}}\rightarrow\mathbb{R}$ being a known $\mathcal{C}^{1}$ function, and $w$ is generated by the following linear exosystem system
\begin{equation}\label{exosystem01}
\begin{split}
  \dot{w} = S_{b}w
\end{split}
\end{equation}
with $w\in\mathbb{R}^{n_{w}}$ and $S_{b}\in \mathbb{R}^{n_{w}\times n_{w}}$. It is assumed that the reference signal is also generated by a linear exosystem as follows:
\begin{equation}\label{exosystem03}
\begin{split}
    \dot{x}_{0}=S_{a}x_{0}.
\end{split}
\end{equation}
where $x_0 \in\mathbb{R}^{r}$ and $S_{a}\in \mathbb{R}^{r \times r}$.
Let $v=\mbox{col}(x_{0},w)$ and $S=\mbox{diag}(S_{a},S_{b})$. Then we can put \eqref{exosystem01} and \eqref{exosystem03} together as follows:
\begin{equation}\label{exosystem1}
\begin{split}
  \dot{v} = Sv.
\end{split}
\end{equation}
The system (\ref{system1}) and the exosystem (\ref{exosystem1}) together can be viewed as a multi-agent system of $(N
+1)$ agents with (\ref{exosystem1}) as the leader and the $N$
subsystems of (\ref{system1}) as $N$ followers. With respect to the plant (\ref{system1}), the exosystem (\ref{exosystem1}), and a given switching signal
$\sigma(t)$, we can define a time-varying digraph
$\bar{\mathcal{G}}_{\sigma(t)}=(\bar{\mathcal{V}},\bar{\mathcal{E}}_{\sigma(t)})$\footnote{See Appendix A for
a summary of graph.} with  $\bar{\mathcal{V}}=\{0,1,\dots,N\}$ and
$\bar{\mathcal{E}}_{\sigma(t)}\subseteq \bar{\mathcal{V}}\times
\bar{\mathcal{V}}$ for all $t\geq0$, where the node $0$ is associated
with the leader system \eqref{exosystem1}  and the node $i$, $i = 1,\dots,N$,
is associated with the $i$th subsystem of system
\eqref{system1}.
For $i=1,\dots,N$, $j=0,1,\dots,N$, and $i\neq j$, $(j,i) \in
\bar{\mathcal{E}}_{\sigma(t)}$ if and only if  $u_i$ can use the
information of the $j$th subsystem for control at time instant $t$.
Let $\bar{\mathcal{A}}_{\sigma(t)} =[\bar{a}_{ij} (t)]\in \mathbb{R}^{(N+1)\times
(N+1)}$ be the weighted adjacency matrix of $\bar{\mathcal{G}}_{\sigma(t)}$.
Let $\bar{\mathcal{N}}_i(t)=\{j,(j,i)\in \bar{\mathcal{E}}_{\sigma(t)}\}$
denote the neighbor set of agent $i$ at time $t$.
Let  $\mathcal{G}_{\sigma(t)}=(\mathcal{V},\mathcal{E}_\sigma(t))$  be the subgraph of $\bar{\mathcal{G}}_{\sigma(t)}$, where $\mathcal{V}=\{1,\cdots,N\}$ and $\mathcal{E}_{\sigma(t)}\subseteq\mathcal{V}\times\mathcal{V}$ is obtained from $\bar{\mathcal{E}}_{\sigma(t)}$ by removing all edges between the node $0$ and the nodes in $\mathcal{V}$.
Clearly, the case where the network topology is fixed can be viewed as a
special case of switching network topology when the switching index
set contains only one element.

Let us  describe our control law as follows.
\begin{equation}\label{ui1}
  \begin{aligned}
 u_{i}&=h_{i}(x_{i},\zeta_{i},x_{j},\zeta_{j})\\
 \dot{\zeta}_{i}&=l_{i}(x_{i},\zeta_{i},x_{j},\zeta_{j},j\in\bar{\mathcal{N}}_{i}(t)),~i=1,\cdots,N\\
  \end{aligned}
\end{equation}
where $h_{i}$ and $l_{i}$  are some nonlinear functions.

A control law of the form \eqref{ui1} is called a distributed control law since $u_{i}$ only depends on the information of its neighbors and itself. Our problem is described as follows.
\begin{Problem}\label{CORPS}
 Given the multi-agent system composed of (\ref{system1}) and (\ref{exosystem1}), and  a switching graph $\bar{\mathcal{G}}_{\sigma(t)}$, design a control law of the form \eqref{ui1},  such that, for any
initial states $x_{i}(0)$, $\zeta_{i}(0)$ and $v(0)$,  the solution of the closed-loop system  exists  for all $t\geq0$, and satisfies
$\lim_{t\rightarrow+\infty}(x_{i}(t)-x_{0}(t))=0$.
\end{Problem}

To solve our problem, we introduce two assumptions as follows.


\begin{Assumption}\label{Ass1}
 All the eigenvalues of $S$ are distinct with zero real parts.
\end{Assumption}
\begin{Remark}\label{RemarkAss1}
Under Assumption \ref{Ass1}, the exosystem \eqref{exosystem1} can generate arbitrarily large constant signals and  multi-tone sinusoidal signals with arbitrarily unknown initial phases and amplitudes  and arbitrarily known frequencies. Since, under Assumption \ref{Ass1},  all the eigenvalues of $S_{a}$ are distinct, the minimal polynomial of $S_{a}$ is equal to the characteristic polynomial of $S_{a}$. Thus, without loss of generality, we can always assume
\begin{equation}\label{Sa1}
\begin{split}
S_{a}=\left[
                                                                                  \begin{array}{cccc}
                                                                                    0 & 1 & \cdots & 0 \\
                                                                                    \vdots & \vdots & \ddots & \vdots \\
                                                                                    0 & 0 & \cdots & 1 \\
                                                                                    \alpha_{1} & \alpha_{2} & \cdots & \alpha_{r} \\
                                                                                  \end{array}
                                                                                \right]
\end{split}
\end{equation}
where $\alpha_{1},\alpha_{2},\cdots,\alpha_{r}$ are some constants. Let $x_{0}=\mbox{col}(x_{10}$, $\cdots,x_{r0})$. Then, we have
\begin{equation}\label{exosystem02}
\begin{split}
    \dot{x}_{s0}& =x_{(s+1)0},~s=1,2\cdots,r-1\\
    \dot{x}_{r0}& = \alpha_{1}x_{10}+\alpha_{2}x_{20}+\cdots+\alpha_{r}x_{r0}.\\
\end{split}
\end{equation}
It is also noted that, under Assumption \ref{Ass1}, given any compact set $\mathbb{V}_{0}$, there exists a compact set $\mathbb{V}$ such that, for any $v(0)\in\mathbb{V}_{0}$, the trajectory $v(t)$ of the exosystem \eqref{exosystem1} remains in $\mathbb{V}$ for all $t\geq0$.
\end{Remark}
\begin{Assumption}\label{Ass3}
There exists a subsequence $\{i_{k}\}$ of $\{i:i=0,1,\cdots\}$ with $t_{i_{k+1}}-t_{i_{k}}<\epsilon$ for some positive $\epsilon$ such that  the union graph $\bigcup_{j=i_{k}}^{i_{k+1}-1}\mathcal{\bar{G}}_{\sigma(t_{j})}$ contains a directed spanning tree with node $0$ as the root.
\end{Assumption}
\begin{Remark}\label{RemarkAss3}
Assumption \ref{Ass3} is  called jointly connected condition in \cite{CH1,LH4,SH5}, which allows the network to be disconnected at any time instant.
\end{Remark}

\section{Adaptive Distributed Observer}\label{Observer}
The key of our approach is to utilize a so-called adaptive distributed observer proposed in \cite{CH1}.
Let us first recall the distributed observer for the leader system of the form (\ref{exosystem1}) as follows \cite{SH5}:
\begin{equation}\label{hatvi3}
\begin{aligned}
\dot{\hat{v}}_{i}&=S\hat{v}_{i}+\mu_{0}\sum_{j=0}^{N}\bar{a}_{ij}(t)(\hat{v}_{j}-\hat{v}_{i}),~i=1,\cdots,N
\end{aligned}
\end{equation}
where $\hat{v}_{0}=v \in \mathbb{R}^{q}$ with $q=r+n_{w}$,  $\hat{v}_{i}\in\mathbb{R}^{q}$  for $i=1,\cdots,N$,  and $\mu_{0}$ is any positive constant.   By Lemma 2 of \cite{SH5}, under Assumptions \ref{Ass1} and \ref{Ass3}, we have $\lim_{t\rightarrow+\infty}(\hat{v}_{i}-v)=0$, $i=1,\cdots,N$. That is why we call \eqref{hatvi3} the distributed observer for \eqref{exosystem1}.

However, a drawback of \eqref{hatvi3} is that the matrix $S$ is used by every follower.
 To overcome this drawback, an adaptive distributed observer for \eqref{exosystem1} was further proposed in \cite{CH1} as follows:
\begin{equation}\label{hatvi4}
\begin{aligned}
\dot{\hat{S}}_{i}&=\mu_{1}\sum_{j=0}^{N}\bar{a}_{ij}(t)(\hat{S}_{j}-\hat{S}_{i})\\
\dot{\hat{v}}_{i}&=\hat{S}_{i}\hat{v}_{i}+\mu_{2}\sum_{j=0}^{N}\bar{a}_{ij}(t)(\hat{v}_{j}-\hat{v}_{i}),~i=1,\cdots,N
\end{aligned}
\end{equation}
where $\hat{v}_{0}=v\in \mathbb{R}^{q}$, $\hat{S}_{0}=S\in \mathbb{R}^{q\times q}$, for $i=1,\cdots,N$, $\hat{v}_{i} \in\mathbb{R}^{q}$,
$\hat{S}_i \in \mathbb{R}^{q \times q} $,  and $\mu_{1}$ and $\mu_{2}$ are any positive constants.

\begin{Remark}
In (\ref{hatvi4}),  the quantity $\hat{v}_{i}$ is to estimate $v$ and the quantity $\hat{S}_{i}$ is to estimate $S$. This is why it is called an adaptive distributed observer. It is noted that $\dot{\hat{S}}_i $  depends on $S$  at time $t$ iff the leader is the neighbor of the $i$th follower at time $t$. Thus, it is
 more practical than the distributed observer proposed in \cite{SH5}
 since  the matrix $S$ is used by every follower in \cite{SH5}.
\end{Remark}

 Let $\tilde{v}_{i}=\hat{v}_{i}-v$ and $\tilde{S}_{i}=\hat{S}_{i}-S$ for $i=0,1,\cdots,N$. Then, for $i=1,\cdots,N$,
 \begin{equation}\label{tildevi3}
\begin{aligned}
\dot{\tilde{S}}_{i}&=\!\mu_{1}\sum_{j=0}^{N}\bar{a}_{ij}(t)(\tilde{S}_{j}-\tilde{S}_{i})\\
\dot{\tilde{v}}_{i}&=\!\tilde{S}_{i}\hat{v}_{i}\!+\!S\tilde{v}_{i}\!+\!\mu_{2}\sum_{j=0}^{N}\bar{a}_{ij}(t)(\tilde{v}_{j}\!-\!\tilde{v}_{i}).
\end{aligned}
\end{equation}
Let $\tilde{v}=\mbox{col}(\tilde{v}_{1},\cdots,\tilde{v}_{N})$, $\hat{v}=\mbox{col}(\hat{v}_{1},\cdots,\hat{v}_{N})$, $\tilde{S}=\mbox{col}(\tilde{S}_{1},\cdots,\tilde{S}_{N})$, and $\tilde{S}_{d}=\mbox{block~diag}\{\tilde{S}_{1},\cdots,\tilde{S}_{N}\}$. Then
\eqref{tildevi3} can be further put into the following compact form
\begin{equation}\label{dottildeSv2}
\begin{split}
&\dot{\tilde{S}}=-\mu_{1}(H_{\sigma(t)}\otimes I_{q})\tilde{S}\\
&\dot{\tilde{v}}=(I_{N}\otimes S-\mu_{2} (H_{\sigma(t)}\otimes I_{q}))\tilde{v}+\tilde{S}_{d}\hat{v}\\
\end{split}
\end{equation}
where $H_{\sigma(t)}=[h_{ij}(t)]_{i,j=1}^{N}$ with $h_{ij}(t)=-\bar{a}_{ij}(t)$ for $i\neq j$ and $h_{ii}(t)=\sum_{j=0}^{N}\bar{a}_{ij}(t)$.

It was shown in Lemmas 1 and 2 of \cite{CH1} that, under Assumptions \ref{Ass1} and \ref{Ass3}
 and the assumption that the subgraph $\mathcal{G}_{\sigma(t)}$ of $\bar{\mathcal{G}}_{\sigma(t)}$ is undirected,
$\lim_{t\rightarrow+\infty}\tilde{S}(t)=0$ exponentially and $\lim_{t\rightarrow+\infty}\tilde{v}(t)=0$ asymptotically.
More recently, the assumption that the subgraph $\mathcal{G}_{\sigma(t)}$ of $\bar{\mathcal{G}}_{\sigma(t)}$ is undirected has been removed in Lemma 2 of \cite{LH4}.
However, to handle the external disturbance $w$, we require
$\lim_{t\rightarrow+\infty}\tilde{v}(t)=0$ exponentially.
For this purpose, we will strengthen  Lemma 2 of \cite{LH4} to the following.
\begin{Lemma}\label{Lemma1}
Under Assumptions \ref{Ass1} and \ref{Ass3},   for any initial conditions $\tilde{S}(0)$ and $\tilde{v}(0)$ and any $\mu_1, \mu_2 > 0$, we have
\begin{enumerate}
  \item $\lim_{t\rightarrow+\infty}\tilde{S}(t)=0$ exponentially;
  \item $\lim_{t\rightarrow+\infty}\tilde{v}(t)=0$ exponentially.
\end{enumerate}
\end{Lemma}
\begin{Proof}
We first note that, Lemma 2 of \cite{LH4} has shown $\lim_{t\rightarrow+\infty}\tilde{S}(t)=0$ exponentially for any $\mu_1 > 0$ by applying Corollary 4 of \cite{SH5}.

Thus, we only need to show $\lim_{t\rightarrow+\infty}\tilde{v}(t)=0$ exponentially.  Note that, in Lemma 2 of \cite{LH4}, we have already shown $\lim_{t\rightarrow+\infty}\tilde{v}(t)=0$ asymptotically. 
Here we further strengthen the result from the asymptotical convergence to the exponential convergence. For convenience, we use the same notation as that in the proof of Lemma 2 of \cite{LH4}.

 Let   $A(t)=(I_{N}\otimes S-\mu_{2} (H_{\sigma(t)}\otimes I_{q}))$ and $F(t)=\tilde{S}_{d}(t) (\mathbf{1}_N \otimes v)$.
Then, the second equation of \eqref{dottildeSv2} can be put to the following form:
 \begin{equation}\label{dottildev2}
\begin{split}
 \dot{\tilde{v}}=A(t)\tilde{v}+\tilde{S}_{d}(t) \tilde{v} +  F (t).
\end{split}
\end{equation}

Note that $H_{\sigma(t)}$ is a piecewise constant matrix with the range of $\sigma$ being $\mathcal{P}=\{1,2,\dots,n_{0}\}$. Thus,
$A(t)$ is bounded over $[0,+\infty)$ and continuous on each time interval $[t_{i},t_{i+1})$ for $i=0,1,2,\cdots$.
 Also note that $v (t)$ is bounded for all time under Assumption \ref{Ass1} and  $\lim_{t\rightarrow+\infty}\tilde{S}(t)=0$ exponentially. Thus,   $F(t)$ is continuous for all $t\geq0$ and $\lim_{t\rightarrow+\infty}F(t)=0$ exponentially. 

By Lemma 2 of \cite{SH5},  under Assumptions \ref{Ass1} and \ref{Ass3},  for any $\mu_{2} > 0$, the origin of the system
 \begin{equation}\label{dottildev1}
\begin{split}
 \dot{\tilde{v}}=A(t)\tilde{v}
\end{split}
\end{equation}
is exponentially stable. Let $\Phi(\tau,t)$ be the state transition matrix of the system (\ref{dottildev1}). Then, we have, for any $t, \tau \geq 0$,
    \begin{equation}\label{dotPhi}
\begin{split}
\frac{\partial}{\partial t}\Phi(\tau,t)=-\Phi(\tau,t)A(t),~\Phi(t,t)=I_{Nq}.
\end{split}
\end{equation}
Since the equilibrium point $\tilde{v}=0$ of \eqref{dottildev1} is exponentially stable, there exist some positive constants $\alpha$ and $\lambda$ such that
    \begin{equation}\label{Phi}
\begin{split}
\|\Phi(\tau,t)\|\leq\alpha e^{-\lambda(\tau-t)}, ~\forall \tau\geq t\geq0.
\end{split}
\end{equation}


Define
$P(t)=\int_{t}^{\infty}\Phi^{T}(\tau,t)Q\Phi(\tau,t)d\tau$,
where  $Q$ is an arbitrarily chosen symmetric positive definite constant matrix.
Then, $P(t)$ is continuous for all $t\geq0$.
Also, it is easy to verify  that there exist some  positive constants $c_{1}$ and $c_{2}$ such that
\begin{equation*}
\begin{split}
 c_{1}\|\tilde{v}\|^{2}\leq \tilde{v}^{T}P(t)\tilde{v}\leq c_{2}\|\tilde{v}\|^{2},
\end{split}
\end{equation*}
which implies that $P(t)$ is positive definite and bounded. Thus there exists a positive constant $c_{3}$ such that $\|P(t)\|\leq c_{3}$ for all $t\geq 0$.

Also, similar to the proof of Theorem 4.12 of \cite{Khalil1}, we have, for any $t\in[t_{i},t_{i+1})$ with $i=0,1,2,\cdots$,
   \begin{equation}\label{dotPt}
\begin{split}
\dot{P}(t)=&\int_{t}^{\infty}\Phi^{T}(\tau,t)Q[\frac{\partial}{\partial t}\Phi(\tau,t)]d\tau+\int_{t}^{\infty}[\frac{\partial}{\partial t}\Phi^{T}(\tau,t)]Q\Phi(\tau,t)d\tau-Q\\
=&-P(t)A(t)-A^{T}(t)P(t)-Q.
\end{split}
\end{equation}

Let $\bar{V}(t,\tilde{v})=\tilde{v}^{T}P(t)\tilde{v}$. For simplicity, we denote $\bar{V}(t,\tilde{v})$ by $\bar{V}(t)$. Clearly, $\bar{V}(t)$ is  positive definite and proper. According to \eqref{dotPt},  for any $t\in[t_{i},t_{i+1})$ with $i=0,1,2,\cdots$, we have
   \begin{equation*}
\begin{split}
\dot{\bar{V}}(t)|_{\eqref{dottildev2}}\!=&\tilde{v}^{T}\dot{P}(t)\tilde{v}+\tilde{v}^{T}(A^{T}(t)P(t)+P(t)A(t))\tilde{v}+  2\tilde{v}^{T}P(t)\tilde{S}_{d}(t) \tilde{v} +  2\tilde{v}^{T}P(t) F (t)  \\
=&\tilde{v}^{T}(-P(t)A(t)-A^{T}(t)P(t)-Q)\tilde{v}+\tilde{v}^{T}(A^{T}(t)P(t)+P(t)A(t))\tilde{v}\\
&+  2\tilde{v}^{T}P(t)\tilde{S}_{d}(t) \tilde{v} +  2\tilde{v}^{T}P(t) F (t)  \\
=&-\!\tilde{v}^{T}Q\tilde{v}\!+\! 2\tilde{v}^{T}P(t)\tilde{S}_{d}(t) \tilde{v} \! +\!  2\tilde{v}^{T}P(t) F (t)  \\
\leq&\!\!-\!\lambda_{\min}(Q)\|\tilde{v}\|^{2}\!\!+\! 2c_3 \|\tilde{S}_d (t)\|  \|\tilde{v}\|^{2} +  \frac{\|P(t)\|^{2}}{\varepsilon}\|\tilde{v}\|^{2}\!\!+\!\varepsilon\|F(t)\|^{2}\\
\leq&-(\lambda_{\min}(Q) - 2c_3 \|\tilde{S}_d (t)\|  -\frac{c_{3}^{2}}{\varepsilon})\|\tilde{v}\|^{2}+\varepsilon\|F(t)\|^{2}.
\end{split}
\end{equation*}
Choose $\varepsilon=\frac{2c_{3}^{2}}{\lambda_{\min}(Q)}$. Then
   \begin{equation*}
\begin{split}
\lambda_{\min}(Q) - 2c_3 \|\tilde{S}_d (t)\|  -\frac{c_{3}^{2}}{\varepsilon}&=\frac{1}{2} \lambda_{\min}(Q) - 2c_3 \|\tilde{S}_d (t)\|.\\
\end{split}
\end{equation*}
Since $\lim_{t\rightarrow+\infty}\tilde{S}_{d}(t)=0$ exponentially,  there  exist some positive integer $l$ and some positive real number $c_{4}$ such that
 \begin{equation*}
\lambda_{\min}(Q) - 2c_3 \|\tilde{S}_d (t)\|  -\frac{c_{3}^{2}}{\varepsilon} > c_{4}>0
\end{equation*}
for all $t \geq t_l$.

Let $\lambda_{1}=\frac{c_{4}}{c_{2}}$. Then,  for any $t\in[t_{i},t_{i+1})$ with $i=l, l+1,l+2,\cdots$, we have
   \begin{equation}\label{dotV2}
\begin{split}
\dot{\bar{V}}(t)|_{\eqref{dottildev2}}
\leq&-c_{4}\|\tilde{v}\|^{2}+\varepsilon\|F(t)\|^{2}\\
\leq&-\lambda_{1}\bar{V}(t)|_{\eqref{dottildev2}}+\varepsilon\|F(t)\|^{2}.\\
\end{split}
\end{equation}
Since $F(t)$ converges to zero exponentially, there exist some positive constants $\gamma_{2}$ and $\lambda_{2}\neq\lambda_{1}$ such that, for any $t\geq t_{l}$,
   \begin{equation}\label{Ft1}
\begin{split}
\varepsilon\|F(t)\|^{2}\leq\gamma_{2}e^{-\lambda_{2}(t-t_{l})}\|F(t_{l})\|^{2}.
\end{split}
\end{equation}
By \eqref{dotV2},  for any $t\geq t_{l}$, we have
   \begin{equation}\label{V1}
\begin{split}
\bar{V}(t)|_{\eqref{dottildev2}}\leq&e^{-\lambda_{1}(t-t_{l})}\bar{V}(t_{l})|_{\eqref{dottildev2}}+\int_{t_{l}}^{t}e^{-\lambda_{1}(t-\tau)}\varepsilon\|F(\tau)\|^{2}d\tau.\\
\end{split}
\end{equation}
According to \eqref{Ft1},
\begin{equation}\label{Ftau1}
\begin{split}
\int_{t_{l}}^{t}e^{-\lambda_{1}(t-\tau)}\varepsilon\|F(\tau)\|^{2}d\tau
&\leq \int_{t_{l}}^{t}e^{-\lambda_{1}(t-\tau)}\gamma_{2}e^{-\lambda_{2}(\tau-t_{l})}\|F(t_{l})\|^{2}d\tau\\
&=\gamma_{2}\|F(t_{l})\|^{2}e^{\lambda_{2}t_{l}}e^{-\lambda_{1}t}\int_{t_{l}}^{t}e^{(\lambda_{1}-\lambda_{2})\tau}d\tau\\
&=\frac{\gamma_{2}\|F(t_{l})\|^{2}}{\lambda_{1}-\lambda_{2}}e^{\lambda_{2}t_{l}}e^{-\lambda_{1}t}(e^{(\lambda_{1}-\lambda_{2})t}-e^{(\lambda_{1}-\lambda_{2})t_{l}})\\
&=\frac{\gamma_{2}\|F(t_{l})\|^{2}}{\lambda_{1}-\lambda_{2}}(e^{-\lambda_{2}(t-t_{l})}-e^{-\lambda_{1}(t-t_{l})}).\\
\end{split}
\end{equation}
 Now let $W(t_{l})=\frac{\gamma_{2}\|F(t_{l})\|^{2}}{|\lambda_{1}-\lambda_{2}|}$ and $\lambda_{0}=\min\{\lambda_{1},\lambda_{2}\}$. Then, according to \eqref{V1} and \eqref{Ftau1}, for any $t\geq t_{l}$,
   \begin{equation}\label{V2}
\begin{split}
\bar{V}(t)|_{\eqref{dottildev2}}
\leq&e^{-\lambda_{1}(t-t_{l})}\bar{V}(t_{l})|_{\eqref{dottildev2}}+(e^{-\lambda_{1}(t-t_{l})}+e^{-\lambda_{2}(t-t_{l})})W(t_{l})\\
\leq&e^{-\lambda_{0}(t-t_{l})}(\bar{V}(t_{l})|_{\eqref{dottildev2}}+2W(t_{l}))
\end{split}
\end{equation}
that is to say, $\lim_{t\rightarrow+\infty}\bar{V}(t)|_{\eqref{dottildev2}}=0$ exponentially. Together with $c_{1}\|\tilde{v}\|^{2}\leq \bar{V}(t)$, we have $\lim_{t\rightarrow+\infty}\tilde{v}(t)|_{\eqref{dottildev2}}=0$ exponentially.
Thus the proof is completed.
\end{Proof}

\begin{Remark}
Lemma \ref{Lemma1}
generalizes Lemma 2 of \cite{LH4} in the sense that we have established that $\lim_{t\rightarrow+\infty}\tilde{v}(t)=0$ exponentially, while, in \cite{LH4}, it was only proved that $\lim_{t\rightarrow+\infty}\tilde{v}(t)=0$ asymptotically. It will be seen that the exponential convergence of $\tilde{v}$ in Lemma \ref{Lemma1} will play a key role in establishing our main result in the next section. Also note that, the main difference between the proof of Lemma \ref{Lemma1} and the proof of Lemma 2 of \cite{LH4} starts from the inequality (\ref{dotV2}).
\end{Remark}

\section{Main Result}\label{MR}
In this section, we will consider the leader-following consensus problem for the system \eqref{system1} subject to jointly connected switching network.
To apply Lemma 3.1 to the leader system (\ref{exosystem1}), where $v = \mbox{col} (x_0, w)$ and $S = \mbox{diag} (S_a, S_b)$, we let
$\hat{v}_{i}=\mbox{col}(\hat{x}_{i},\hat{w}_{i})\in\mathbb{R}^{q}$ with $\hat{x}_{i}\in\mathbb{R}^{r}$ and $\hat{w}_{i}\in\mathbb{R}^{n_{w}}$ denoting the estimations of $x_{0}$ and $w$ respectively, $\hat{S}_{i}=\mbox{diag}(\hat{S}_{ai},\hat{S}_{bi})\in\mathbb{R}^{q\times q}$ with
\begin{equation*}
\begin{aligned}
\hat{S}_{ai}=\left[
                                                                                  \begin{array}{cccc}
                                                                                    0 & 1 & \cdots & 0 \\
                                                                                    \vdots & \vdots & \ddots & \vdots \\
                                                                                    0 & 0 & \cdots & 1 \\
                                                                                    \hat{\alpha}_{1i} & \hat{\alpha}_{2i} & \cdots & \hat{\alpha}_{ri} \\
                                                                                  \end{array}
                                                                                \right]\in\mathbb{R}^{r\times r}
\end{aligned}
\end{equation*}
 and $\hat{S}_{bi}\in\mathbb{R}^{n_{w}\times n_{w}}$
denoting the estimations of $S_{a}$ and $S_{b}$ respectively. Finally, let $\hat{x}_{i}=\mbox{col}(\hat{x}_{1i},\cdots,\hat{x}_{ri})$.
Then, we can decompose the adaptive distributed observer \eqref{hatvi4} into two parts as follows: for $s=1,\cdots,r-1,$ and $i=1,\cdots,N$,
\begin{equation}\label{hatxi2}
\begin{aligned}
\dot{\hat{S}}_{ai}&=\mu_{1}\sum_{j=0}^{N}\bar{a}_{ij}(t)(\hat{S}_{aj}-\hat{S}_{ai})\\
\dot{\hat{x}}_{si}&=\!\hat{x}_{(s+1)i}\!+\!\mu_{2}\sum_{j=0}^{N}\bar{a}_{ij}(t)(\hat{x}_{sj}\!-\!\hat{x}_{si})\\
\dot{\hat{x}}_{ri}&=\!\sum_{s=1}^{r}\hat{\alpha}_{si}\hat{x}_{si}\!+\!\mu_{2}\sum_{j=0}^{N}\bar{a}_{ij}(t)(\hat{x}_{rj}\!-\!\hat{x}_{ri})\\
\end{aligned}
\end{equation}
and
\begin{equation}\label{hatwi2}
\begin{aligned}
\dot{\hat{S}}_{bi}&=\mu_{1}\sum_{j=0}^{N}\bar{a}_{ij}(t)(\hat{S}_{bj}-\hat{S}_{bi})\\
\dot{\hat{w}}_{i}&=\hat{S}_{bi}\hat{w}_{i}+\mu_{2}\sum_{j=0}^{N}\bar{a}_{ij}(t)(\hat{w}_{j}-\hat{w}_{i}).\\
\end{aligned}
\end{equation}


\begin{Remark}\label{Remark4.2}
By Lemma \ref{Lemma1}, 
under Assumptions \ref{Ass1} and \ref{Ass3}, for $i=1,\cdots,N$, we have $\lim_{t\rightarrow+\infty}\tilde{S}_{i}(t)=0$ and $\lim_{t\rightarrow+\infty}\tilde{v}_{i}(t)=0$, which implies that
\begin{equation}\label{hatxik3}
\begin{split}
\lim_{t\rightarrow+\infty}(\hat{x}_{ki}(t)-x_{k0}(t))=0,~k=1,\cdots,r.
\end{split}
\end{equation}
On the other hand, since $\hat{v}_{i}$ and $\bar{a}_{ij}(t)$ are bounded for all $t\geq0$, we obtain that  $\lim_{t\rightarrow+\infty}\dot{\tilde{S}}_{i}(t)=0$ and $\lim_{t\rightarrow+\infty}\dot{\tilde{v}}_{i}(t)=0$. Let $\tilde{x}_{dsi}=\mu_{2}\sum_{j=0}^{N}\bar{a}_{ij}(t)(\hat{x}_{sj}\!-\!\hat{x}_{si})$ 
 for $s=1,\cdots,r$ and $i=1,\cdots,N$.  Then it is easy to verify that
\begin{equation}\label{hatxik4}
\begin{split}
\lim_{t\rightarrow+\infty}\tilde{x}_{dsi}(t)=0.
\end{split}
\end{equation}

\end{Remark}

Next, we will develop a distributed control law. For this purpose, we let
\begin{equation}\label{pri1}
\begin{split}
p_{ri}=&~\hat{x}_{ri}\!-\!\beta_{1}(x_{(r\!-\!1)i}\!-\!\hat{x}_{(r\!-\!1)i})\!-\!\cdots-\beta_{r\!-\!1}(x_{1i}\!-\!\hat{x}_{1i})\\
\end{split}
\end{equation}
where $\beta_{1},\cdots,\beta_{r-1}$ are some positive constants such that the polynomial $\lambda^{r-1}+\beta_{1}\lambda^{r-2}+\cdots+\beta_{r-2}\lambda+\beta_{r-1}=0$ is stable. Then
\begin{equation}\label{dotpri1}
\begin{split}
\dot{p}_{ri}=&\dot{\hat{x}}_{ri}\!-\!\beta_{1}(\dot{x}_{(r\!-\!1)i}\!-\!\dot{\hat{x}}_{(r\!-\!1)i})\!-\!\cdots-\beta_{r\!-\!1}(\dot{x}_{1i}\!-\!\dot{\hat{x}}_{1i})\\
=&\dot{\hat{x}}_{ri}\!-\!\beta_{1}(x_{ri}\!-\!\dot{\hat{x}}_{(r\!-\!1)i})\!-\!\cdots-\beta_{r\!-\!1}(x_{2i}\!-\!\dot{\hat{x}}_{1i}).\\
\end{split}
\end{equation}
Let
\begin{equation}\label{si2}
\begin{split}
s_{i}&=x_{ri}-p_{ri}.
\end{split}
\end{equation}
Now we are ready to present our control law as follows:
\begin{equation}\label{ui3}
\begin{split}
u_{i}&=-f_{i}^{T}(x_{i},t)\hat{\theta}_{i}-d_{i}(\hat{w}_{i})-k_{i}s_{i}+\dot{p}_{ri}\\
\dot{\hat{\theta}}_{i}&=\Lambda_{i}^{-1}f_{i}(x_{i},t)s_{i}\\
\dot{\hat{S}}_{i}&=\mu_{1}\sum_{j=0}^{N}\bar{a}_{ij}(t)(\hat{S}_{j}-\hat{S}_{i})\\
\dot{\hat{v}}_{i}&=\hat{S}_{i}\hat{v}_{i}+\mu_{2}\sum_{j=0}^{N}\bar{a}_{ij}(t)(\hat{v}_{j}-\hat{v}_{i}),~i=1,\cdots,N
\end{split}
\end{equation}
where $k_{i}$ is some positive constant and $\Lambda_{i}\in\mathbb{R}^{m\times m}$ is some symmetric positive definite matrix.

The closed-loop system composed of \eqref{system1} and \eqref{ui3} is as follows:
 \begin{equation}\label{system3}
\begin{split}
    \dot{x}_{si}& =x_{(s+1)i},~s=1,2\cdots,r-1\\
    \dot{x}_{ri}& = -f_{i}^{T}(x_{i},t)\tilde{\theta}_{i}+\tilde{d}_{i}(\tilde{w}_{i},w)-k_{i}s_{i}+\dot{p}_{ri}\\
    \dot{\hat{\theta}}_{i}&=\Lambda_{i}^{-1}f_{i}(x_{i},t)s_{i}\\
    \dot{\tilde{S}}_{i}&=\mu_{1}\sum_{j=0}^{N}\bar{a}_{ij}(t)(\tilde{S}_{j}-\tilde{S}_{i})\\
\dot{\tilde{v}}_{i}&=\tilde{S}_{i}\hat{v}_{i}+S\tilde{v}_{i}+\mu_{2}\sum_{j=0}^{N}\bar{a}_{ij}(t)(\tilde{v}_{j}-\tilde{v}_{i})\\
\end{split}
\end{equation}
where $\tilde{\theta}_{i}=\hat{\theta}_{i}-\theta_{i}$, $\tilde{w}_{i}=\hat{w}_{i}-w$ and $\tilde{d}_{i}(\tilde{w}_{i},w)=d_{i}(w)-d_{i}(\hat{w}_{i})=d_{i}(w)-d_{i}(\tilde{w}_{i}+w)$.

\begin{Remark}
The construction of the control law (\ref{ui3}) is based on the certainty equivalence principle. Suppose the state  of the leader system is available by the control law of every follower.
Then, instead of (\ref{pri1}), we can define $p_{ri}$ as follows:
\begin{equation*}
\begin{split}
p_{ri}\!=\!{x}_{r0}\!-\!\beta_{1}(x_{(r\!-\!1)i}\!-\!{x}_{(r\!-\!1)0})\!-\!\cdots-\beta_{r\!-\!1}(x_{1i}\!-\!{x}_{10})\\
\end{split}
\end{equation*}
Then, by the standard adaptive control method as can be found, say, in \cite{slotine}, we can show that the following so-called decentralized control law
\begin{equation}\label{ui30}
\begin{split}
u_{i}&=-f_{i}^{T}(x_{i},t)\hat{\theta}_{i}-d_{i}({w})-k_{i}s_{i}+\dot{p}_{ri}\\
\dot{\hat{\theta}}_{i}&=\Lambda_{i}^{-1}f_{i}(x_{i},t)s_{i},~i=1,\cdots,N\\
\end{split}
\end{equation}
solves Problem \ref{CORPS}. However, the control law (\ref{ui30}) is impractical because  $w$ is in general not available for control and it is uninteresting to assume
the state $x_0$ of the leader system is available by every follower. To overcome this difficulty, we replace the leader's state in (\ref{ui30}) by its estimation $\hat{v}_{i}$ generated by
the adaptive distributed observer, thus leading to our distributed control law (\ref{ui3}). Since the distributed control law (\ref{ui3}) is more complicated than the decentralized
control law (\ref{ui30}), we need to provide a much more sophisticated stability analysis for the closed-loop system (\ref{system3}).
\end{Remark}

Now we give our result as follows.
\begin{Theorem}\label{Theorem2}
Under Assumptions \ref{Ass1} and \ref{Ass3}, the leader-following consensus problem for the multi-agent system composed of \eqref{system1} and \eqref{exosystem1} is solvable by the distributed control law \eqref{ui3}.
\end{Theorem}
\begin{Proof}
Note that the closed-loop system \eqref{system3} is piecewise continuous. Thus the traditional Barbalat's Lemma can not be used to analyze the stability of the  closed-loop system \eqref{system3}, and we need to resort to the general Barbalat's Lemma proposed in \cite{SH3}. In order to apply the general Barbalat's Lemma (i.e. Lemma 1 of \cite{SH3}), we need to find a scalar continuous function $V:[0,+\infty)\rightarrow \mathbb{R}$ such that
\begin{enumerate}
  \item $\lim_{t\rightarrow+\infty} V(t)$ exists;
  \item $V(t)$ is twice differentiable on each interval $[t_{i},t_{i+1})$;
  \item $\ddot{V}(t)$ is bounded over $[0,+\infty)$.
\end{enumerate}
Then $\dot{V}(t)\rightarrow 0$ as $t\rightarrow +\infty$. Thus the main challenge here is to find a function $V(t)$ with the above three conditions satisfied.

 Since $\tilde{d}_{i}(0,w)=0$ for all $w$ and $\tilde{d}_{i}$ is $C^1$, by Lemma 11.1 of  \cite{ChenHuang2015}, there exists some smooth function $\bar{d}_{i}(\tilde{w}_{i},w)\geq0$, such that, for all $w\in\mathbb{R}^{n_{w}}$, 
 \begin{equation}\label{tildew2}
\begin{split}
\tilde{d}_{i}^{2}(\tilde{w}_{i},w)\leq \bar{d}_{i}(\tilde{w}_{i},w)\|\tilde{w}_{i}\|^{2}.
\end{split}
\end{equation}
Let 
 \begin{equation}\label{V3}
\begin{split}
V=\frac{1}{2}\sum_{i=1}^{N}(s_{i}^{2}+\tilde{\theta}_{i}^{T}\Lambda_{i}\tilde{\theta}_{i}).
\end{split}
\end{equation}
Then the time derivative of $V$ along the trajectory of the closed-loop system \eqref{system3} is given by
 \begin{equation}\label{dotV3}
\begin{split}
\dot{V}=&\sum_{i=1}^{N}(s_{i}\dot{s}_{i}+\tilde{\theta}_{i}^{T}\Lambda_{i}\dot{\tilde{\theta}}_{i})\\
=&\sum_{i=1}^{N}(s_{i}(\dot{x}_{ri}-\dot{p}_{ri})+\tilde{\theta}_{i}^{T}\Lambda_{i}\dot{\hat{\theta}}_{i})\\
=& \sum_{i=1}^{N}\bigg(s_{i}(-f_{i}^{T}(x_{i},t)\tilde{\theta}_{i}+\tilde{d}_{i}(\tilde{w}_{i},w)-k_{i}s_{i})+\tilde{\theta}_{i}^{T}f_{i}(x_{i},t)s_{i}\bigg)\\
=&\sum_{i=1}^{N}\bigg(s_{i}\tilde{d}_{i}(\tilde{w}_{i},w)-k_{i}s_{i}^{2}\bigg)\\
\leq&\sum_{i=1}^{N}\bigg(\frac{1}{4}s_{i}^{2}+\tilde{d}_{i}^{2}(\tilde{w}_{i},w)-k_{i}s_{i}^{2}\bigg)\\
\leq&\sum_{i=1}^{N}\bigg(-(k_{i}-\frac{1}{4})s_{i}^{2}+\bar{d}_{i}(\tilde{w}_{i},w)\|\tilde{w}_{i}\|^{2}\bigg)\\
\end{split}
\end{equation}
Choosing $k_{i}\geq\frac{5}{4}$ gives
 \begin{equation}\label{dotV4}
\begin{split}
\dot{V}\leq&\sum_{i=1}^{N}\bigg(-s_{i}^{2}+\bar{d}_{i}(\tilde{w}_{i},w)\|\tilde{w}_{i}\|^{2}\bigg)\\
\end{split}
\end{equation}
 and thus
\begin{equation}\label{V4}
\begin{split}
V(t)=&\int_{0}^{t}\dot{V}(\tau)d\tau+c_{0}\\
\leq&\int_{0}^{t}\sum_{i=1}^{N}\bar{d}_{i}(\tilde{w}_{i},w)\|\tilde{w}_{i}\|^{2}d\tau +c_{0}
\end{split}
\end{equation}
where $c_{0}$ is some constant. The existence of the limit $\lim_{t\rightarrow+\infty}V(t)$ is shown in Appendix B.
Moreover, by Lemma \ref{Lemma1}, under Assumptions \ref{Ass1} and \ref{Ass3}, $\lim_{t\rightarrow+\infty}\tilde{v}(t)=0$ exponentially,  which implies $\lim_{t\rightarrow+\infty}\tilde{w}_{i}(t)=0$ exponentially. Thus,  $V(t)$ is bounded for all $t\geq0$. Thus, for $i=1,\cdots,N$, $s_{i}$ and $\tilde{\theta}_{i}$ are bounded for all $t\geq0$. By \eqref{dotpri1} and \eqref{si2}, $s_{i}(t)$ is differentiable  on each interval $[t_{k},t_{k+1})$ for all $k\geq0$ and so is $\dot{V}(t)$. By \eqref{pri1} and \eqref{si2}, we have
\begin{equation}\label{xri3}
\begin{split}
&x_{ri}+\beta_{1}x_{(r-1)i}+\cdots+\beta_{r-1}x_{1i} \\
&=s_{i}+\hat{x}_{ri}+\beta_{1}\hat{x}_{(r-1)i}+\cdots+\beta_{r-1}\hat{x}_{1i}
\end{split}
\end{equation}
which is equivalent to
\begin{equation}\label{dotx1i3}
\begin{split}
&x_{1i}^{(r-1)}+\beta_{1}x_{1i}^{(r-2)}+\cdots+\beta_{r-1}x_{1i} \\
&=s_{i}+\hat{x}_{ri}+\beta_{1}\hat{x}_{(r-1)i}+\cdots+\beta_{r-1}\hat{x}_{1i}.
\end{split}
\end{equation}
Since both $s_{i}$ and $\hat{x}_{i}$ are  bounded,  \eqref{dotx1i3} can be viewed as a stable $(r-1)$th order linear system in $x_{1i}$ with a bounded input, and thus $x_{1i}^{(r-1)},x_{1i}^{(r-2)}\cdots,x_{1i}$ are all bounded, that is to say, $x_{i}$ is bounded. Therefore, from \eqref{pri1} and \eqref{dotpri1},  $p_{ri}$ and $\dot{p}_{ri}$ are both bounded. From the second equation of \eqref{system3},  $\dot{x}_{ri}$ is bounded. Thus $\dot{s}_{i}=\dot{x}_{ri}-\dot{p}_{ri}$ is also bounded. Note that
 \begin{equation}\label{ddotV2}
\begin{split}
\ddot{V}=&\sum_{i=1}^{N}\bigg(s_{i}(\frac{\partial \tilde{d}_{i}(\tilde{w}_{i},w)}{\partial \tilde {w}_{i}}\dot{\tilde{w}}_{i}+\frac{\partial \tilde{d}_{i}(\tilde{w}_{i},w)}{\partial w}\dot{w})+\dot{s}_{i}\tilde{d}_{i}(\tilde{w}_{i},w)-2k_{i}s_{i}\dot{s}_{i}\bigg).
\end{split}
\end{equation}
Since $s_{i}$, $\dot{s}_{i}$, $w$, $\dot{w}$, $\tilde{w}_{i}$ and $\dot{\tilde{w}}_{i}$ are all bounded, there exists a positive number $\gamma$ such that 
 \begin{equation}\label{ddotV3}
\begin{split}
\sup_{t_{k}\leq t<t_{k+1},k=0,1,\cdots}|\ddot{V}(t)|\leq \gamma.
\end{split}
\end{equation}
It follows from the generalized Barbalat's Lemma as noted at the beginning of this proof, we have $\lim_{t\rightarrow+\infty}\dot{V}(t)=0$.  Thus, from \eqref{dotV4}, we further have $\lim_{t\rightarrow+\infty}s_{i}(t)=0$ for $i=1,\cdots,N$.

Next, we let $z_{1i}=x_{1i}-\hat{x}_{1i}$, $z_{2i}=x_{2i}-\hat{x}_{2i}$, $\cdots$, $z_{(r-1)i}=x_{(r-1)i}-\hat{x}_{(r-1)i}$. Then, from \eqref{hatxi2}, \eqref{system3} and \eqref{xri3}, we have
 \begin{equation}\label{system4}
\begin{split}
\dot{z}_{si}=&z_{(s+1)i}-\tilde{x}_{dsi},~s=1,\cdots,r-2\\
\dot{z}_{(r-1)i}=&-\beta_{r-1}z_{1i}-\cdots-\beta_{1}z_{(r-1)i}+s_{i}-\tilde{x}_{d(r-1)i}.\\
\end{split}
\end{equation}
Further, let $z_{i}=\mbox{col}(z_{1i},\cdots,z_{(r-1)i})$ and
 \begin{equation*}
\begin{split}
A=\left[\!\!\!
                    \begin{array}{cccc}
                      0 & 1 & \cdots & 0 \\
                      \vdots & \vdots & \ddots & \vdots \\
                      0 & 0 & \cdots & 1 \\
                      -\beta_{r-1} & -\beta_{r-2} & \cdots & -\beta_{1} \\
                    \end{array}
                  \!\!\!\right]\!,~\bar{u}_{i}=\left[\!\!
                                         \begin{array}{c}
                                           -\tilde{x}_{d1i} \\
                                           \vdots \\
                                           -\tilde{x}_{d(r-2)i} \\
                                           s_{i}\!-\!\tilde{x}_{d(r-1)i} \\
                                         \end{array}
                                       \!\!\right].
\end{split}
\end{equation*}
Then \eqref{system4} can be put into the following form
 \begin{equation}\label{system5}
\begin{split}
\dot{z}_{i}&=Az_{i}+\bar{u}_{i}.
\end{split}
\end{equation}
Since $A$ is Hurwitz, and $\lim_{t\rightarrow+\infty}\bar{u}_{i}(t)=0$ by \eqref{hatxik4}, we obtain that $\lim_{t\rightarrow+\infty}z_{i}(t)=0$ i.e. $\lim_{t\rightarrow+\infty}(x_{si}(t)-\hat{x}_{si}(t))=0$ for $s=1,\cdots,r-1$ and $i=1,\cdots,N$. By \eqref{xri3}, we have
$ \lim_{t\rightarrow+\infty}(x_{ri}(t)-\hat{x}_{ri}(t)=\lim_{t\rightarrow+\infty}(s_{i}-\beta_{1}(x_{(r-1)i}-\hat{x}_{(r-1)i})-\cdots-\beta_{r-1}(x_{1i}-\hat{x}_{1i}))=0$
for $i=1,\cdots,N$.

Together with \eqref{hatxik3}, we have
%
%
%
 \begin{equation}\label{dotxsi2}
\begin{split}
\lim_{t\rightarrow+\infty}(x_{si}(t)-x_{0i}(t))=0,~s=1,\cdots,r.
\end{split}
\end{equation}
Thus our proof is completed.
\end{Proof}
\begin{Remark}
Our proof relies critically on the fact that the signal  $\tilde{w}_i (t)$  converges to zero exponentially, which is established in Lemma 3.1. If the convergence of $\tilde{w}_i (t)$ is only asymptotic but not exponential, then we cannot guarantee that $V (t)$ is bounded over $t \geq 0$.
\end{Remark}

\begin{Remark}
As mentioned in the introduction, problems similar to ours have also been studied via the adaptive control technique in \cite{Das1,LH3,YuH1,ZhangH1}. It is interesting to make some comparisons.  First, references  \cite{Das1,YuH1}  only considered the first-order nonlinear systems and reference \cite{LH3} only considered the second-order nonlinear systems while here we study the higher-order nonlinear systems. Second, in \cite{YuH1}, the systems do not contain external disturbance and the overall control laws are not distributed in the sense that the control law of each follower has to depend on the information of the leader.  Third, references \cite{Das1,ZhangH1} employed  neutral networks to approximate certain unknown nonlinear functions and thus the control laws can only make the tracking errors uniformly ultimately bounded,  and the results are not global in the sense that they are only valid for initial conditions in the prescribed compact set. Finally, references \cite{Das1,LH3,ZhangH1} only considered the static network case which can be viewed as a special case of jointly connected switching network case.
\end{Remark}

\section{An Example}\label{Example}

Consider the leader-following consensus problem for a group of Vol del Pol systems as follows:
\begin{equation}\label{VDP1}
\begin{split}
  \dot{x}_{1i}& =x_{2i},~~ i=1,2,3,4\\
    \dot{x}_{2i}& = -\theta_{1i}x_{1i}+\theta_{2i}x_{2i}(1-x_{1i}^{2})+d_{i}(w) +u_{i}\\
\end{split}
\end{equation}
where
  \begin{equation*}
\begin{split}
&x_{i}=\mbox{col}(x_{1i},x_{2i})\in\mathbb{R}^{2},~~w=[w_{1}, w_{2}]^{T},~~d_{1}(w)=w_{1}^{2}w_{2}^{2}\\
&d_{2}(w)=w_{1}w_{2}^{3},~~d_{3}(w)=w_{1}^{3}+w_{1}w_{2},~~d_{4}(w)=w_{2}^{4}.\\
\end{split}
\end{equation*}
System \eqref{VDP1} is in the form \eqref{system1} with
  \begin{equation*}
\begin{split}
f_{i}(x_{i},t)=
\left[
  \begin{array}{c}
    -x_{1i} \\
    x_{2i}(1-x_{1i}^{2})\\
  \end{array}
\right],~~~
                                                                                         \theta_{i}=\left[
                                                                                                      \begin{array}{c}
                                                                                                         \theta_{1i}\\
                                                                                                        \theta_{2i} \\
                                                                                                      \end{array}
                                                                                                    \right].
\end{split}
\end{equation*}

 The exosystem is in the form (\ref{exosystem1}) with
\begin{equation*}
\begin{split}
&S_{a}=\left[
    \begin{array}{cc}
      0 & 1 \\
      -1 & 0 \\
    \end{array}
  \right],~S_{b}=\left[
    \begin{array}{cc}
      0 & 0.5 \\
      -0.5 & 0 \\
    \end{array}
  \right].
\end{split}
\end{equation*}
It can be seen that Assumption \ref{Ass1} is satisfied.

\begin{figure}[H]
  \centering
  \subfigure[$\bar{\mathcal{G}}_{1}$]{
    \includegraphics[width=0.85in]{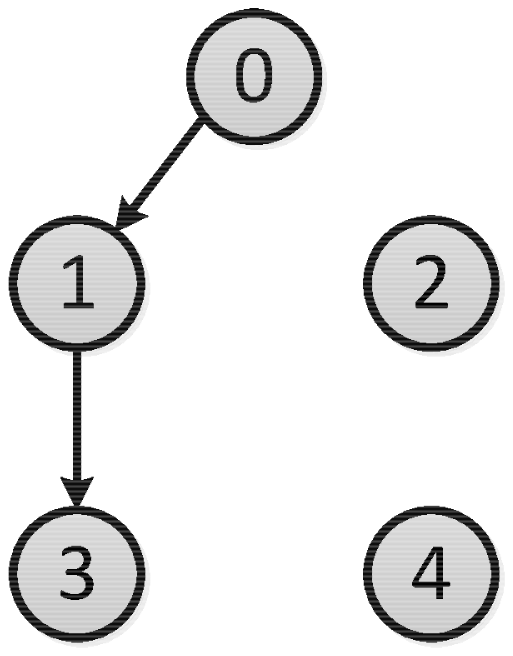}}
      \hspace{0.4in}
  \subfigure[$\bar{\mathcal{G}}_{2}$]{
    \includegraphics[width=0.85in]{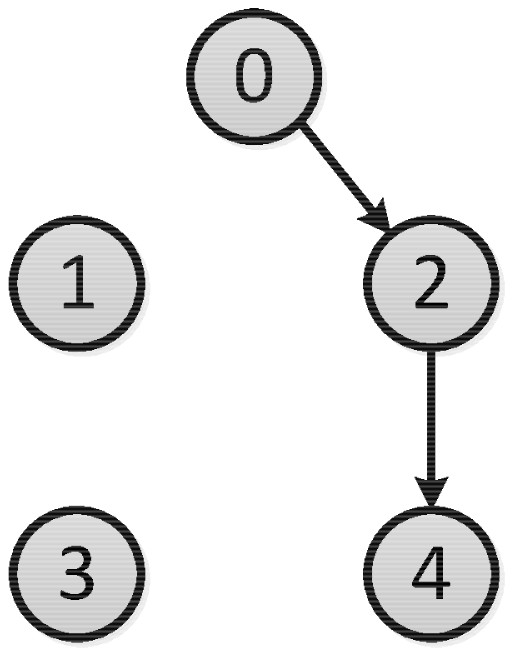}}
          \hspace{0.4in}
  \subfigure[$\bar{\mathcal{G}}_{3}$]{
    \includegraphics[width=0.85in]{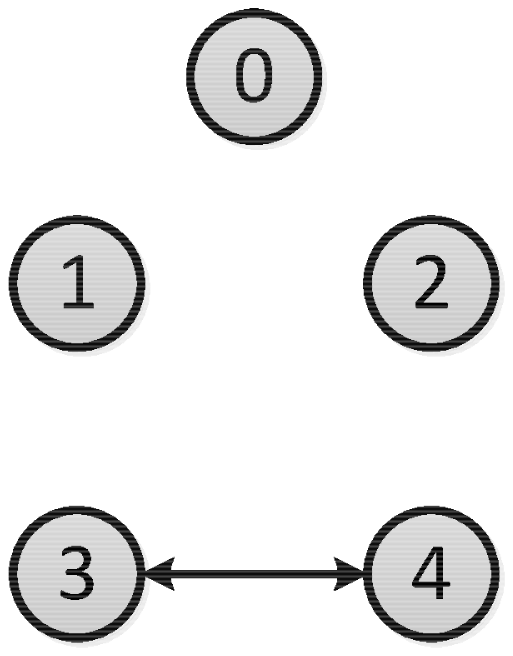}}
     \hspace{0.4in}
  \subfigure[$\bar{\mathcal{G}}_{4}$]{
    \includegraphics[width=0.85in]{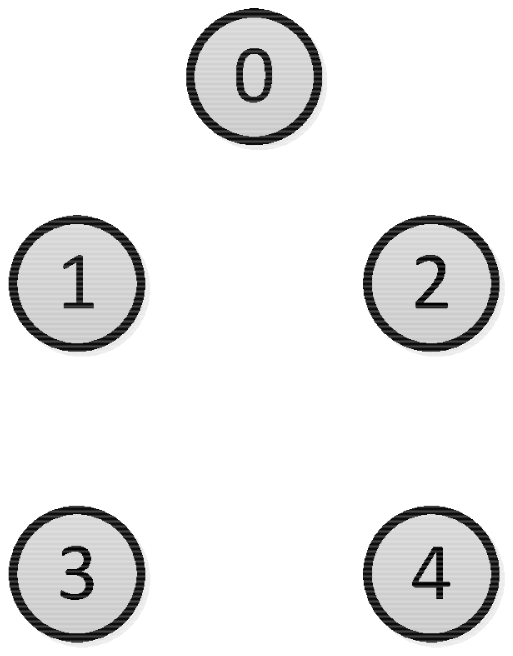}}
  \caption{Switching topology $\bar{\mathcal{G}}_{\sigma(t)}$ with $\mathcal{P}=\{1,2,3,4\}$}
  \label{switching} 
\end{figure}

The communication graph $\bar{\mathcal{G}}_{\sigma(t)}$ is
dictated by the following switching signal:
\begin{equation}\label{sigmat2}
\begin{split}
\sigma(t)= \left\{
  \begin{array}{ll}
    1, & \hbox{if\ $sT_{0}\leq t< (s+\frac{1}{4})T_{0}$} \\
    2, & \hbox{if\ $(s+\frac{1}{4})T_{0}\leq t< (s+\frac{1}{2})T_{0}$} \\
    3, & \hbox{if\ $(s+\frac{1}{2})T_{0}\leq t< (s+\frac{3}{4})T_{0}$}\\
    4, & \hbox{if\ $(s+\frac{3}{4})T_{0}\leq t< (s+1)T_{0}$}\\
  \end{array}
\right.
 \end{split}
\end{equation}
where $s=0,1,2,\cdots$. The four digraphs $\bar{\mathcal{G}}_{i}$, $i=1,2,3,4$, are described by Figure \ref{switching} where the node $0$ is associated with the leader and the other nodes are associated with the followers.
It can be verified that Assumption \ref{Ass3} is satisfied even though the four digraphs $\bar{\mathcal{G}}_{i}$ are all disconnected.

By Theorem \ref{Theorem2}, we can design a distributed control law of the form \eqref{ui3} with $\mu_{1}=3$, $\mu_{2}=12$, $\beta_{1}=1$  and $k_{i}=3$ for $i=1,2,3,4$.

\begin{figure}[H]
\centering
\includegraphics[scale=0.55]{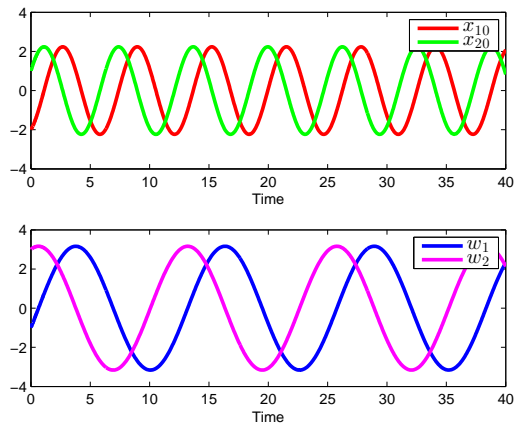}
\caption{States of leader system: $v=\mbox{col}(x_{0},w)$.} \label{leader1}
\end{figure}
\begin{figure}[H]
\centering
\includegraphics[scale=0.55]{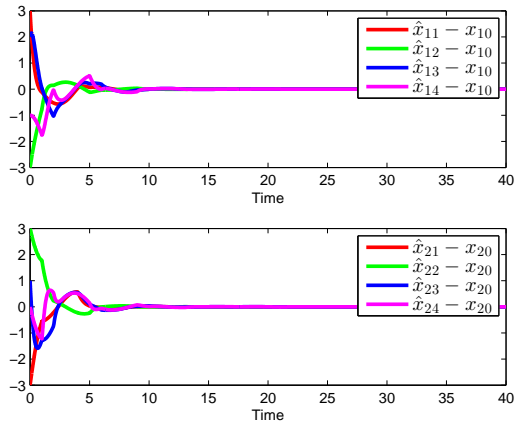}
\caption{Estimation errors: $\hat{x}_{i}-x_{0}$.} \label{estimation1}
\end{figure}
\begin{figure}[H]
\centering
\includegraphics[scale=0.55]{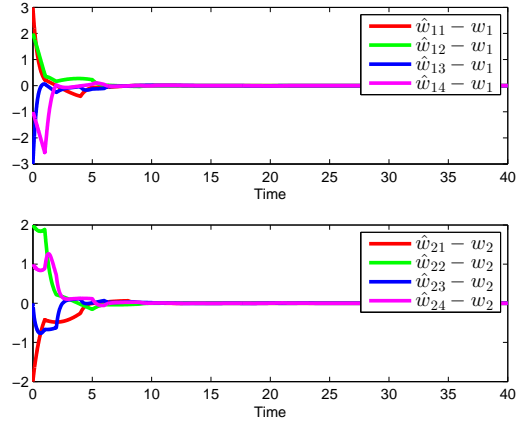}
\caption{Estimation errors: $\hat{w}_{i}-w$.} \label{estimation2}
\end{figure}
\begin{figure}[H]
\centering
\includegraphics[scale=0.55]{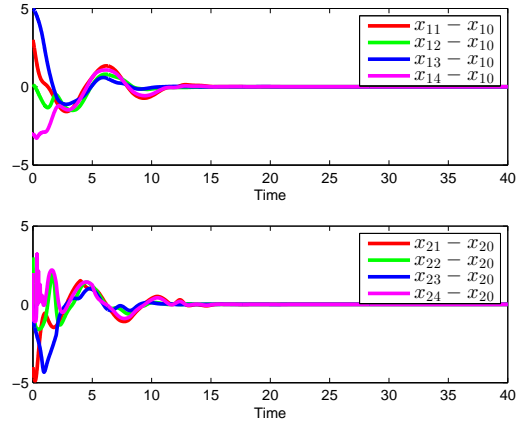}
\caption{Tracking errors: $x_{i}-x_{0}$.} \label{error1}
\end{figure}

Simulation is performed with 
 $\theta_{1}=[4,5]^{T}$, $\theta_{2}=[3,1]^{T}$, $\theta_{3}=[2,5]^{T}$, $\theta_{4}=[5,3]^{T}$,
  and the following initial conditions:
\begin{equation*}
  \begin{aligned}
&x_{1}(0)=[1,-4]^{T},~x_{2}(0)=[-2,3]^{T},~x_{3}(0)=[3,1]^{T}\\
&x_{4}(0)=[-5,2]^{T},~v(0)=[-2,1,-1,3]^{T}\\
 &\hat{v}_{1}(0)=[1, -2,  2, 1]^{T},~  \hat{v}_{2}(0)= [-5, 4, 1, 5]^{T} \\
 &\hat{v}_{3}(0)=   [0, 2, -4, 3]^{T},~  \hat{v}_{4}(0)= [-3, 1, -2, 4]^{T}\\
 &\hat{\theta}_{i}(0)=0_{2\times1},~\hat{S}_{bi}(0)=0_{2\times 2}\\
&\hat{S}_{ai}(0)=\left[
                   \begin{array}{cc}
                     0 & 1 \\
                     0 & 0 \\
                   \end{array}
                 \right]
,~i=1,2,3,4.\\
  \end{aligned}
\end{equation*}

Figure \ref{leader1} shows the states of the leader system. Figure \ref{estimation1} and \ref{estimation2} show the estimation errors between the states of the adaptive distributed observer of each subsystem and the states of the leader system, which approach zero as time tends to infinity. Figure \ref{error1} shows that the states of all followers approach the states of the leader asymptotically.
All these simulations confirm that our control law is effective in solving our problem even though the communication network is switching and disconnected at every time constant.

\section{Conclusion}\label{Conclusion}
In this paper, we have studied the leader-following consensus problem for a class of higher-order nonlinear multi-agent systems subject to both constant parameter uncertainties and external disturbances under jointly connected switching networks. By combining the adaptive control technique and the established technical lemma on the adaptive distributed observer, our problem has been solved by the designed distributed state feedback control law.


\appendix
\subsection{Digraph}
A digraph
${\mathcal{G}}=({\mathcal{V}},{\mathcal{E}})$ consists of a finite set of nodes  ${\mathcal{V}}=\{1\cdots,N\}$ and an edge set  ${\mathcal{E}}\subseteq
{\mathcal{V}}\times{\mathcal{V}}$.
An edge of ${\mathcal{E}}$ from node $i$ to node $j$ is denoted by $(i,j)$,
 where nodes $i$ and $j$ are called the parent node and the child node of each other, and node $i$ is called a neighbor of node $j$.
  Define
$\mathcal{{N}}_{i}=\{j|(j,i)\in{\mathcal{E}} \}$, which is called the neighbor set of node $i$.
The edge $(i, j)$ is called undirected if $(i, j)\in \mathcal{{E}}$ implies $(j, i)\in \mathcal{{E}}$.
 The digraph ${\mathcal{G}}$ is called
undirected if every edge in $\mathcal{{E}}$ is undirected.
If the digraph ${\mathcal {G}} $ contains a sequence of edges of the form $(i_{1},i_{2}),(i_{2},i_{3}),\cdots,(i_{k},i_{k+1})$, then the set $\{(i_{1},i_{2}),(i_{2},i_{3}),\cdots,(i_{k},i_{k+1})\}$ is called a directed path of ${\mathcal {G}} $ from node $i_{1}$ to node $i_{k+1}$ and node $i_{k+1}$ is said to be reachable from node $i_{1}$.
A directed tree is a digraph where every node has exactly one parent node except for one node
called the root, from which every other node is reachable.
 A digraph ${\mathcal {G}}_{s}=({\mathcal{V}}_{s},{\mathcal{E}}_{s})$ is called a subgraph of the digraph ${\mathcal{G}}=({\mathcal{V}},{\mathcal{E}})$ if ${\mathcal{V}}_{s}\subseteq{\mathcal{V}}$ and ${\mathcal{E}}_{s}\subseteq{\mathcal{E}}\cap({\mathcal{V}}_{s}\times{\mathcal{V}}_{s})$.
 A subgraph ${\mathcal {G}}_{s}=({\mathcal{V}}_{s},{\mathcal{E}}_{s})$ of the diagraph ${\mathcal{G}}=({\mathcal{V}},{\mathcal{E}})$ is called a directed spanning tree of ${\mathcal {G}}$ if ${\mathcal {G}}_{s}$ is a directed tree and ${\mathcal{V}}_{s}={\mathcal{V}}$.
 Given a set of $n_{0}$ digraphs $\{{\mathcal{G}}_{i}=({\mathcal{V}},{\mathcal{E}}_{i}),i=1,\cdots,n_{0}\}$, the digraph ${\mathcal{G}}=({\mathcal{V}},{\mathcal{E}})$ where ${\mathcal{E}}=\bigcup_{i=1}^{n_{0}}{\mathcal{E}}_{i}$ is called the union of digraphs ${\mathcal{G}}_{i}$, denoted by ${\mathcal{G}}=\bigcup_{i=1}^{n_{0}}{\mathcal{G}}_{i}$.

The weighted adjacency matrix of the digraph ${\mathcal{G}}$ is a nonnegative matrix
${\mathcal{A}}=[{a}_{ij}]\in \mathbb{R}^{N\times
N}$ where ${a}_{ii}=0$ and ${a}_{ij}>0\Leftrightarrow
(j,i)\in\mathcal{{E}}$, $i,j=1,\cdots,N$.
 On the other hand, given a matrix ${\mathcal{A}}=[{a}_{ij}]\in \mathbb{R}^{N\times N}$ satisfying ${a}_{ii}=0$ and ${a}_{ij}\geq0$ for $i,j=0,1,\cdots,N$, we can always define a digraph ${\mathcal{G}}$ such that ${\mathcal{A}}$ is the weighted adjacency matrix of the digraph ${\mathcal{G}}$. We call ${\mathcal{G}}$ the digraph of ${\mathcal{A}}$.


Given a piecewise
constant switching signal $\sigma:
[0,+\infty)$ $\rightarrow \mathcal{P}=\{1,2,\dots,n_{0}\}$, and a set of
$n_{0}$ graphs $\mathcal{{G}}_i=(\mathcal{{V}},\mathcal{{E}}_i)$, $i=1,\dots,n_{0}$  with the corresponding weighted adjacency
matrices being denoted by $\mathcal{{A}}_i$, $i = 1,\cdots, n_{0}$,  we call a time-varying graph $\mathcal{{G}}_{\sigma
(t)}=(\mathcal{{V}},\mathcal{{E}}_{\sigma (t)})$ a switching graph, and denote the weighted adjacency
matrix of $\mathcal{{G}}_{\sigma
(t)}$ by $\mathcal{{A}}_{\sigma
(t)}$.

\subsection{Existence of the limit $\lim_{t\rightarrow+\infty}V(t)$}
\begin{Proof}
First, let
 \begin{equation*}
\begin{split}
\tilde{z}=\mbox{vec}(\tilde{S})=\left[
                        \begin{array}{c}
                          \tilde{S}_{c1} \\
                          \vdots \\
                          \tilde{S}_{cq} \\
                        \end{array}
                      \right]
\end{split}
\end{equation*}
where $\tilde{S}_{ci}$, $i=1,2,\cdots,q$, is the $i$th column of $\tilde{S}$.
Then the first equation of \eqref{dottildeSv2} can be  put into the following form:
 \begin{equation}\label{dottildez1}
\begin{split}
&\dot{\tilde{z}}=-\mu_{1}(I_{q}\otimes H_{\sigma(t)}\otimes I_{q})\tilde{z}.\\
\end{split}
\end{equation}
By Lemma \ref{Lemma1}, $\lim_{t\rightarrow+\infty}\tilde{z}(t)=0$ exponentially. Then, similar to the construction of the Lyapunov function for \eqref{dottildev1} in the proof of Lemma \ref{Lemma1}, we can also construct a continuous and piecewise differentiable  quadratic Lyapunov function $V_{1}(t,\tilde{z})$ for \eqref{dottildez1} such that
 \begin{equation}\label{V1a}
\begin{split}
l_{1}\|\tilde{z}\|^{2}\leq V_{1}(t,\tilde{z})\leq l_{2}\|\tilde{z}\|^{2}
\end{split}
\end{equation}
 \begin{equation}\label{dotV1a}
\begin{split}
\frac{\partial V_{1}}{\partial t}+\frac{\partial V_{1}}{\partial \tilde{z}} (-\mu_{1}(I_{q}\otimes H_{\sigma(t)}\otimes I_{q})\tilde{z})\leq - l_{3}\|\tilde{z}\|^{2}
\end{split}
\end{equation}
for some positive constants $l_{1},l_{2}$, and $l_{3}$, and  all $t\in[t_{k},t_{k+1})$ with $k=0,1,2,\cdots$.

Second, for convenience, repeat \eqref{dottildev2} as follows.
 \begin{equation}\label{dottildev3}
\begin{split}
 &\dot{\tilde{v}}=A(t)\tilde{v}+\tilde{S}_{d}(t) \tilde{v} + F(t). 
\end{split}
\end{equation}
Let $V_{2}(t,\tilde{v}) $ be the same as the function $\bar{V} (t,\tilde{v})$ in  the proof of Lemma \ref{Lemma1}. Then, as shown in Lemma \ref{Lemma1}, $V_{2}(t,\tilde{v})$ is such that
 \begin{equation}\label{V2a}
\begin{split}
c_{1}\|\tilde{v}\|^{2}\leq V_{2}(t,\tilde{v})\leq c_{2}\|\tilde{v}\|^{2}
\end{split}
\end{equation}
for some positive constants $c_{1},c_{2}$ and all $t\geq0$, and
 \begin{equation}\label{dotV2a}
\begin{split}
\frac{\partial V_{2}}{\partial t}\!+\!\frac{\partial V_{2}}{\partial \tilde{v}} (A(t)\tilde{v}\!+\!\tilde{S}_{d}(t) \tilde{v} \!+\! F(t))\leq& -\! c_{4}\|\tilde{v}\|^{2}+\varepsilon\|F(t)\|^{2}
\end{split}
\end{equation}
for some positive constants $c_{4},\varepsilon$, and  all $t\in[t_{k},t_{k+1})$ with $k=l,l+1,\cdots$, where $l$ is some positive integer.

Since $\|F(t)\|=\|\tilde{S}_{d}(t) (\mathbf{1}_N \otimes v(t))\|$ and, under Assumption \ref{Ass1}, $v(t)$ is bounded for all $t\geq0$, there  exists a constant $l_{4}$ such that $\varepsilon\|\tilde{S}_{d}(\mathbf{1}_N \otimes v)\|^{2}\leq l_{4}\|\mbox{vec}(\tilde{S})\|^{2}=l_{4}\|\tilde{z}\|^{2}$ for all $t\geq0$. Thus we have
 \begin{equation}\label{dotV2b}
\begin{split}
\frac{\partial V_{2}}{\partial t}\!+\!\frac{\partial V_{2}}{\partial \tilde{v}} (A(t)\tilde{v}\!+\!\tilde{S}_{d}(t) \tilde{v} \!+\! F(t))\leq& -\! c_{4}\|\tilde{v}\|^{2}+l_{4}\|\tilde{z}\|^{2}
\end{split}
\end{equation}
for all $t\in[t_{k},t_{k+1})$ with $k=l,l+1,\cdots$.

For convenience, we further let $\tilde{X}=\mbox{col}(\tilde{z},\tilde{v})$ and $f(t,\tilde{X})=\mbox{col}(-\mu_{1}(I_{q}\otimes H_{\sigma(t)}\otimes I_{q})\tilde{z},A(t)\tilde{v}+\tilde{S}_{d}(t) \tilde{v} +  F (t))$. Then \eqref{dottildez1} and \eqref{dottildev3} can be put together into the following compact form:
 \begin{equation}\label{dottildeX1}
\begin{split}
 \dot{\tilde{X}}=f(t,\tilde{X}).
\end{split}
\end{equation}
Let $V_{3}(t,\tilde{X})=l_{5}V_{1}(t,\tilde{z})+V_{2}(t,\tilde{v})$ with $l_{5}\geq\frac{l_{4}+c_{4}}{l_{3}}$ being a positive constant. Then, by \eqref{V1a} and \eqref{V2a}, for all $t\geq0$, we have
 \begin{equation}\label{V3a}
\begin{split}
m_{1}\|\tilde{X}\|^{2}\leq V_{3}(t,\tilde{X})\leq m_{2}\|\tilde{X}\|^{2}
\end{split}
\end{equation}
where $m_{1}=\min\{l_{1}l_{5},c_{1}\}$ and $m_{2}=\max\{l_{2}l_{5},c_{2}\}$.
 Moreover, according to \eqref{dotV1a} and \eqref{dotV2b}, we have
 \begin{equation}\label{dotV3a}
\begin{split}
\frac{\partial V_{3}}{\partial t}\!\!+\!\!\frac{\partial V_{3}}{\partial \tilde{X}} f(t,\tilde{X})\!\leq& \!-\!l_{3}l_{5}\|\tilde{z}\|^{2}\!-\! c_{4}\|\tilde{v}\|^{2}\!+\!l_{4}\|\tilde{z}\|^{2}\\
\leq&-c_{4}\|\tilde{z}\|^{2}- c_{4}\|\tilde{v}\|^{2}\\
=&-c_{4}\|\tilde{X}\|^{2}
\end{split}
\end{equation}
for all $t\in[t_{k},t_{k+1})$ with $k=l,l+1,\cdots$.

Note that, under Assumption \ref{Ass1}, by Remark \ref{RemarkAss1}, there exist a compact subset $\mathbb{W}\subset\mathbb{R}^{n_{w}}$ such that $w(t)\in \mathbb{W}$ for all $t\geq0$.
Then there exists some smooth positive  function $\hat{d}(\|\tilde{X}\|) $ such that, for all $w\in\mathbb{W}$,
 \begin{equation}\label{tildew3}
\begin{split}
 \sum_{i=1}^{N}\bar{d}_{i}(\tilde{w}_{i},w)\|\tilde{w}_{i}\|^{2}\leq \hat{d}(\|\tilde{X}\|)\|\tilde{X}\|^{2}.
\end{split}
\end{equation}

By Lemma 11.2 of \cite{ChenHuang2015}, there exists a smooth non-decreasing function $\rho:[0,+\infty)\rightarrow[0,+\infty)$ satisfying $\rho(s)>0$ for all $s>0$, such that
 \begin{equation}\label{rho1}
\begin{split}
\rho(m_{1}\|\tilde{X}\|^{2})\|\tilde{X}\|^{2}\geq\frac{1}{c_{4}}\hat{d}(\|\tilde{X}\|)\|\tilde{X}\|^{2}.
\end{split}
\end{equation}
Together with \eqref{V3a}, we have
\begin{equation}\label{rho2}
\begin{split}
\rho(V_{3})\|\tilde{X}\|^{2}\geq\rho(m_{1}\|\tilde{X}\|^{2})\|\tilde{X}\|^{2}\geq\frac{1}{c_{4}}\hat{d}(\|\tilde{X}\|)\|\tilde{X}\|^{2}.
\end{split}
\end{equation}
Let $\bar{V}_{3}=\int_{0}^{V_{3}}\rho(s)ds$. Clearly, $\bar{V}_{3}$ is also a proper and positive definite function.  Then, according to \eqref{dotV3a} and \eqref{rho2}, we have
\begin{equation}\label{dotbV1}
\begin{split}
\dot{\bar{V}}_{3}|_{\eqref{dottildeX1}}&= \rho(V_{3})\bigg(\frac{\partial V_{3}}{\partial t}+\frac{\partial V_{3}}{\partial \tilde{X}} f(t,\tilde{X})\bigg)\\
&\leq  - c_{4}\rho(V_{3})\|\tilde{X}\|^{2}\\
&\leq  - c_{4}\rho(m_{1}\|\tilde{X}\|^{2})\|\tilde{X}\|^{2}\\
&\leq-\hat{d}(\|\tilde{X}\|)\|\tilde{X}\|^{2}
\end{split}
\end{equation}
for all $t\in[t_{k},t_{k+1})$ with $k=l,l+1,\cdots$.
Since $\bar{V}_{3}$ is lower bounded and $\dot{\bar{V}}_{3}|_{\eqref{dottildeX1}}$ is non-positive, we can conclude that $\lim_{t\rightarrow+\infty}\bar{V}_{3}(t)$ exists.

On the other hand, recall from \eqref{V3} that $V=\frac{1}{2}\sum_{i=1}^{N}(s_{i}^{2}+\tilde{\theta}_{i}^{T}\Lambda_{i}\tilde{\theta}_{i})$.
Then \eqref{dotV4} and \eqref{tildew3} imply
 \begin{equation}\label{dotV5}
\begin{split}
 &\dot{V}|_{\eqref{system3}}\leq-\sum_{i=1}^{N}s_{i}^{2}+\hat{d}(\|\tilde{X}\|)\|\tilde{X}\|^{2}.\\
\end{split}
\end{equation}

Let $U=V+\bar{V}_{3}$. Clearly, $U$ is also a proper and positive definite function. Then, by  \eqref{dotbV1} and \eqref{dotV5}, we have
\begin{equation}\label{dotU1}
\begin{split}
\dot{U}|_{\eqref{system3},\eqref{dottildeX1}}\leq& -\sum_{i=1}^{N}s_{i}^{2}+\hat{d}(\|\tilde{X}\|)\|\tilde{X}\|^{2}-\hat{d}(\|\tilde{X}\|)\|\tilde{X}\|^{2}\\
=&-\sum_{i=1}^{N}s_{i}^{2}\\
\end{split}
\end{equation}
for all $t\in[t_{k},t_{k+1})$ with $k=l,l+1,\cdots$.
Since $U$ is lower bounded and $\dot{U}$ is non-positive, we  conclude that $\lim_{t\rightarrow+\infty}U(t)$ exists.

Thus $\lim_{t\rightarrow+\infty}V(t)=\lim_{t\rightarrow+\infty}U(t)-\lim_{t\rightarrow+\infty}\bar{V}_{3}(t)$ also exists.
\end{Proof}

\end{document}